\documentclass[preprint,12pt]{elsarticle}

\usepackage{amssymb}
\usepackage{amsmath}
\usepackage{amsthm}

\usepackage{lineno}

\usepackage{hyperref}
\usepackage{algorithm,algorithmic}

\newtheorem{theorem}{Theorem}[section]
\newtheorem{lemma}[theorem]{Lemma}

\newtheorem{identity}{Identity}[section]
\newtheorem{definition}{Definition}[section]

\newtheorem{example}[definition]{Example}

\journal{Advances in Applied Mathematics}

\begin{document}

\begin{frontmatter}


\title{$q$-Binomial Identities Finder}
\author[label3]{Hao Zhong\corref{cor1}}
\ead{11435011@zju.edu.cn}
\cortext[cor1]{Corresponding author}
\affiliation[label3]{organization={College of Computer Science and Cyber Security, Chengdu University},
            addressline={1 Dongsan Road},
            city={Chengdu},
            postcode={610059},
            state={Sichuan},
            country={China}}




\author[label3]{Leqi Zhao} 


\begin{abstract}
This paper presents a symbolic computation method for automatically transforming $q$-hypergeometric identities to $q$-binomial identities. Through this method, many previously proven $q$-binomial identities, including $q$-Saalsch{\"u}tz’s formula and $q$-Suranyi’s formula, are re-fund, and numerous new ones are discovered. Moreover,
the generation of the identities is accompanied by the corresponding proofs.
During the transformation process, different ranges of variable values and various combinations of $q$-Pochhammer symbols yield different identities. The algorithm maps variable constraints to positive elements in an ordered vector space and employs a backtracking method to provide the feasible variable constraints and $q$-binomial coefficient combinations for each step.
\end{abstract}

%

%


\end{frontmatter}


\section{Introduction}\label{sec:introduction}

The study of Gaussian binomial coefficients, also known as $q$-binomial coefficients, along with their numerous related identities, is valuable in various mathematical and applied contexts due to their broad range of applications and deep theoretical implications. These fields include number theory \cite{andrews1986q}, Lie algebra \cite{feinsilver1990lie}, quantum physics \cite{kupershmidt2010mathematics, jacob2012quantum} and even privacy protection \cite{bewong2019privacy}. The main goal of this paper is to introduce an idea that is both ``old'' and ``new'' for generating more $q$-binomial identities. We refer to it as ``old'' because it is rooted in the traditional method of deriving $q$-binomial identities by rearranging $q$-Pochhammer symbols in $q$-hypergeometric identities, such as the $q$-Gauss summation formula \cite{andrews1974applications, gould1997some}. However, it is also ``new'' as it has successfully automated the generation of identities with the assistance of computers. Before delving into more details, let us start with some common notations in $q$ series.

Throughout this paper, $q$ always represents a complex number such that $|q|<1$. For a positive integer $n$, the $q$-Pochhammer symbol is defined as 
$$(a;q)_{n}=\prod_{k=0}^{n-1}(1-aq^{k}).$$
The $q$-Pochhammer symbol can be naturally extended to zero, infinite and negative products as follows.

$$(a;q)_{0}=1, \quad (a;q)_{\infty}=\prod_{k\geq0}(1-aq^{k}), \quad (a;q)_{-n}=\frac{(a;q)_{\infty}}{(aq^{-n};q)_{\infty}}.$$

Let $m$ and $n$ be two integers. The $q$-binomial coefficient is defined as the following combination of $q$-Pochhammer symbols.
\begin{equation}\label{def:gbc}
\binom{m+n}{n}_{q}=
\begin{cases}
    \frac{(q;q)_{m+n}}{(q;q)_{m}(q;q)_{n}}  &\text{if } m\text{ and }n\geq 0,\\
    0 &\text{otherwise.}
\end{cases}
\end{equation}
The $q$-binomial coefficients are set to zero if either $m$ or $n$ is negative, which makes the coefficients not always match $\frac{(q;q)_{m+n}}{(q;q)_{m}(q;q)_{n}}$. Consequently, $q$-binomial identities hold true only within a fixed range of variables \cite{schilling2000generalization}.

The first significant $q$-binomial identity is the $q$-analog of Vandermonde convolution, which connects $q$-binomial coefficients to the summation of products of $q$-binomial coefficients.
\begin{identity}[$q$-Vandermonde convolution, \cite{gould1997some, slater1966generalized}]\label{id:vandermonde}
For non-negative integers $m$ and $n$,
\begin{equation}
\sum_{r = 0}^{\min(n,k)}q^{r(m-k+r)}{\binom{m}{k-r}}_{q}{\binom{n}{r}}_{q}={\binom{m+n}{k}}_{q}.
\end{equation}
\end{identity}
Identity \ref{id:vandermonde} is proved by Heine (or $q$-Gauss) Summation Theorem \cite{gould1997some} or combinatorially proved by counting finite vector subspaces of certain objects \cite{slater1966generalized}.

The most studied $q$-binomial identity is the following Saalsch{\"u}tz’s formula of $q$-binomial form.
\begin{identity}[Gould, \cite{gould1972new, carlitz1974remark, andrews1975identities, goulden1985bijective2}]\label{id:saaschiitz}
For non-negative integers $a$ and $b$,
\begin{multline}
\sum_{r = 0}^{\min(a,b)}q^{(a-r)(b-r)}{\binom{x+y+r}{r}}_{q}{\binom{x+a-b}{a-r}}_{q}{\binom{y+b-a}{b-r}}_{q}\\={\binom{x+a}{a}}_{q}{\binom{y+b}{b}}_{q}.   
\end{multline}
\end{identity}

Identity \ref{id:saaschiitz} is proved by applying $q$-Vandermonde convolution several times \cite{gould1972new}, or combinatorially proved by the enumeration of ordered pairs of subsets of $\{1,2,3,\ldots,v\}$ \cite{andrews1975identities}, or combinatorially proved by determining an explicit bijection between sets of integer partitions \cite{goulden1985bijective2}.

Besides Identity \ref{id:saaschiitz}, there exist plenty of $q$-binomial identities with ``three [$q$-binomials] in the sum side'' and ``two [$q$-binomials] in the product side'', such as the following identity.
\begin{identity} [$q$-analog of Suranyi’s formula, \cite{pic1965combinatorial}]\label{id:suranyi}
For non-negative integers $k$ and $n$,
\begin{equation}
\sum_{r = 0}^{m}q^{r^2}{\binom{k}{r}}_{q}{\binom{n}{r}}_{q}{\binom{x+k+n-r}{k+n}}_{q}={\binom{x+k}{k}}_{q}{\binom{x+n}{n}}_{q}.
\end{equation}
where $m=\min(k,n,x+k+n)$
\end{identity}

Identity \ref{id:suranyi} is proved by extensions of the ordinary Vandermonde theorem \cite{pic1965combinatorial}.

From the proof techniques mentioned above, we can observe that, the basic hypergeometric series, also known as $q$-hypergeometric series, has been a fertile ground for discovering $q$-binomial identities. 
The most common way to generate Gaussian binomial identities is to make some unpredictable substitutions in $q$-hypergeometric identities and then to simplify them by carefully matching and reorganizing terms on both sides of the identity until the expressions can be rewritten entirely in terms of $q$-binomials, with all extraneous $q$-Pochhammer symbols eliminated.
Here, $q$-hypergeometric function is defined as follows.
$$
{}_r\phi_s \left( \begin{array}{c}
a_1, a_2, \dots, a_r \\
b_1, b_2, \dots, b_s
\end{array}; q, z \right)
= \sum_{n=0}^{\infty} \frac{(a_1, a_2, \dots, a_r; q)_n}{(q, b_1, b_2, \dots, b_s; q)_n} \left( (-1)^n q^{\binom{n}{2}} \right)^{s-r+1} z^n
$$
where $\left(a_{1},a_{2},\cdots,a_{r};q\right)_{n}=\prod_{j=1}^{r}(a_{j};q)_{n}$.
This approach, albeit straightforward, is particularly laborious. Determining combinations of variable substitutions often relies on chance, and the simplification process requires both keen insight and continuous experimentation. 
The only way to avoid relying on luck is to test all possible combinations, which leads to significantly increased manual calculations. 

To address the complexity of manual computations, computer algebra has been widely applied in the study of partition and $q$-series identities \cite{andrews1986q}. For example, Kanade and Russell \cite{kanade2015identityfinder} proposed six challenging partition identity conjectures of Rogers–Ramanujan type using symbolic computation. Then Kur{\c{s}}ung{\"o}z \cite{kurcsungoz2019andrews}, Kanade and Russell \cite{kanade2018staircases} found the generating functions of these partitions in a combinatorial way. Chern and Li \cite{chern2020linked} later used computer algebra and the concept of linked partition ideals \cite{andrews1998theory} to rediscover these six generating function identities. For more computer application in $q$-series, please refer to \cite{zeilberger1987q, frye2019automatic, WANG2024109743}.
This motivates us to use computational tools to overcome the tedious work in transforming $q$-hypergeometric identities to $q$-binomial identities. Consequently, we introduce a Python-based method in this paper. The source code is publicly available at \url{https://github.com/Z798971901/q_binomial_identities_finder} and has also been archived on Zenodo at \url{https://doi.org/10.5281/zenodo.15867344}.  
The method can be summarized into three main modules. 
\begin{enumerate}
    \item \textbf{Constraint Conditions Generator}: Ensuring each step of transformation is valid by restricting the range of variables.
    \item \textbf{$q$-Binomial Coefficients Generator}: Converting each $q$-Pochhammer symbol in the original $q$-hypergeometric identity into a square bracket form of positive integers, and then transform them to $q$-binomial coefficients.
    \item \textbf{Backtracking Framework}: Trying all possible constraint conditions of variables and as many as possible transformations from square brackets to $q$-binomial coefficients.
\end{enumerate}

Although not all $q$-hypergeometric identities can be transformed into $q$-binomial identities, when such a transformation is feasible, our identity finder can efficiently generate $q$-binomial identities. This method has two main applications.

First, our finder can generate $q$-binomial identities. In addition to known identities such as Identity \ref{id:saaschiitz} and Identity \ref{id:suranyi}, we can discover new ones by selecting less-studied $q$-hypergeometric identities as input. For example, feeding Bailey–Daum $q$-Kummer sum into our finder results in the following identity.
\begin{identity}\label{id:kummer}
For two positive integers $m$ and $n$,
\begin{multline}
\sum_{r = 0}^{n}q^{\binom{r+1}{2}}{\binom{2m+r-1}{r}}_{q}{\binom{2m+2n}{n-r}}_{q}{\binom{m+n}{n}}_{q^2}\\=\left(1-q^{n+1}\right){\binom{2m+2n}{2n+1}}_{q}{\binom{2n+1}{n}}_{q}.
\end{multline}
\end{identity}

Second, following each step of our algorithm can provide proofs for the output $q$-binomial identities. This is because our finder essentially simulates a rigorous human proof, with the only difference being that it assembles $q$-binomial coefficients faster and more comprehensively. This will be demonstrated in the subsequent sections.


The remainder of this paper is organized as follows. Section \ref{section:preliminaries} introduces the necessary preliminaries related to square brackets and positive spans; Section \ref{section:methods} provides a detailed description of our method; Section \ref{section:results} lists some of the identities discovered by our generator, which lead to Identity \ref{id:saaschiitz}, Identity \ref{id:suranyi} and Identity \ref{id:kummer}; and finally, Section \ref{section:further works} concludes this paper with a discussion of potential future works.

\section{Preliminaries}\label{section:preliminaries}
\subsection{Square Brackets}
Square brackets serve as a unified language throughout our $q$-identities generator. To clarify their role and application, we will first define square brackets in detail, explaining how they function as an essential tool for translating $q$-Pochhammer symbols to $q$-binomial coefficients.

\begin{definition}[Square Bracket of $q$-Pochhammer symbol]\label{def:sym}
    Let $k$ and $n$ be two non-negative integers, and let $m$ be an integer. We define the square bracket notation as follows.
    $$[m;k]_{n} := (q^{km};q^{k})_{n}$$
    For the sake of simplicity, we denote $[m;1]_{n}$ by $[m]_{n}$, $[m;k]_{\infty}$ by $[m;k]$ and $[m;1]_{\infty}$ by $[m]$, respectively. Additionally, we call $[m]_{n}$ the square bracket with subscript $n$ and $[m]$ the pure square bracket.
\end{definition}

To represent more $q$-Pochhammer symbols using square brackets, we present the following lemmas.

\begin{lemma}\label{lem:extent def}
Follow the notations in Definition \ref{def:sym}. Then
\begin{eqnarray}
(-q^{km};q^{k})_{n}=\frac{[m;2k]_{n}}{[m;k]_{n}},\\
\prod_{r=1}^{k-1}(q^{km+r};q^{k})_{n}=\frac{[km;1]_{kn}}{[m;k]_{n}}.\label{eq:general sym}
\end{eqnarray}
\end{lemma}
This lemma directly follows from Definition \ref{def:sym}, so we omit the proof.

\begin{lemma}\label{lem:2inf}
Let $m$ be an integer and let $n$ be a non-negative integer. Then
\begin{equation}
[m]_{n}=
\begin{cases}
    (-1)^{n}q^{\left(mn+\binom{n}{2}\right)}\frac{[-m-n+1]}{[-m+1]}  &\text{if } m \leq -n,\\
    0 &\text{if } -n < m \leq 0,\\
    \frac{[m]}{[m+n]}  &\text{if } m > 0.
\end{cases}
\end{equation}
Moreover, $[m]_{n}$ can be represented using square brackets of positive integers, or simply equals to zero.
\end{lemma}

An important fact is that the Gaussian binomial symbols can be also represented using square brackets of positive integers, or simply equals to zero.
\begin{lemma}\label{lem:qbinom}
Let $m$ and $n$ be two integers. Then
\begin{equation}
{\binom{m+n}{m}}_{q}=
\begin{cases}
    \frac{[m+1,n+1]}{[m+n+1,1]} &\text{if } m\text{ and }n\geq 0,\\
    0 &\text{otherwise}.
\end{cases}
\end{equation}
\end{lemma}

Now our main purpose can be formally expressed by using the square bracket symbols of positive integers as to transform the following $q$-hypergeomtric identity
\begin{equation}\label{eq:q-hyper}
\sum_{n \geq 0}\alpha(n)\prod_{k=1}^{\infty}\frac{[a_{1},\ldots,a_{r_{k}};k]}{[b_{1},\ldots,b_{s_{k}};k]}=\prod_{k=1}^{\infty}\frac{[c_{1},\ldots,c_{t_{k}};k]}{[d_{1},\ldots,d_{u_{k}};k]}
\end{equation}
to as many as possible $q$-binomial identities
\begin{equation}\label{eq:q-binom}
\sum_{n \geq 0}\beta(n)\prod_{k=1}^{\infty}{\binom{A_{1}}{B_{1}}_{q^{k}}\cdots{\binom{A_{v_{k}}}{B_{v_{k}}}}}_{q^{k}}=\prod_{k=1}^{\infty}{\binom{C_{1}}{D_{1}}_{q^{k}}\cdots{\binom{C_{w_{k}}}{D_{w_{k}}}}}_{q^{k}}.
\end{equation}

\subsection{Positivity}
Determining whether an item inside the square bracket is positive or negative might not pose a significant challenge to a human, but it can be a complex task for a computer at times. This is particularly true when the constraints of variables come into play, making it difficult to determine the polarity of an expression. Simply storing positive or negative results in memory also prevents computers from detecting contradictions in subsequent calculations. In Section \ref{section:methods}, we will address this issue. For now, let us lay the groundwork by introducing the definition of ordered vector space from \cite{hausner1952ordered}.

\begin{definition} [Ordered vector space]\label{def:order}
Given a vector space $V$ over the real numbers $\mathbb{R}$, and a strict partial order $<$ on the set $V$. We say $<$ is compatible with $V$, and $(V,<)$ is an ordered vector space if for any $\textbf{u}$ and $\textbf{v}$ in $V$,
\begin{enumerate}[(1)]
    \item $\textbf{u}<\textbf{0}$, $r$ is real and positive, implies $r\textbf{u}<\textbf{0}$;
    \item $\textbf{u}<\textbf{0}$, $\textbf{v}<\textbf{0}$, implies $\textbf{u}+\textbf{v}<\textbf{0}$;
    \item $\textbf{u}<\textbf{v}$ if and only if $\textbf{u}-\textbf{v}<\textbf{0}$.
\end{enumerate}
\end{definition}
Consequently, $\textbf{0} < \textbf{u}$, or simply denoted as $\textbf{u} > \textbf{0}$ if and only if $-\textbf{u} < \textbf{0}$. Moreover, we denote $\leq$ (respectively, $\geq$) as the associated non-strict partial order relation of $<$ (respectively, $>$). For any subset $S$ of $V$, let $S_{+}:=\{\textbf{u}\in S:\textbf{u}\ge \textbf{0}\}$. Then $V_{+}$ is a convex cone, with $\textbf{0}$ being its vertex as well as its infimum.

Having defined the order concept of a vector space $V$, we now turn our attention to the specific case where $V=\mathbb{R}^{n}$. Each $n$-variable affine expression can be associated with a vector in $V=\mathbb{R}^{n+1}$ where the coefficients of the affine expression correspond one-to-one with the components of the vector. Specifically, the affine expression $a_{0}+a_{1}x_{1}+\cdots+a_{n}x_{n}$ can be associated with the vector $(a_{0},a_{1},\ldots,a_{n})$ in $V=\mathbb{R}^{n+1}$. Obviously, this one-to-one correspondence is an isomorphism with respect to addition. Thus, the compatible order $<$ can also be inherited. Once 
$\mathbb{R}^{n+1}_{+}$ is determined, we can obtain a polytope that defines the range of values for $x$'s. To ensure the existence of integer solutions, we first add an infimum to $\mathbb{Z}^{n}_{+}\backslash\{\textbf{0}\}$. Specifically, We extend the definition of $<$ by requiring $\textbf{u}-\textbf{1}\geq0$ for any nonzero vector $u$ in $\mathbb{Z}^{n}_{+}$, where $\textbf{1}=(1,0,0,\ldots,0)$. Then, we slightly modify the definition of positive span and frame as follows.

\begin{definition}[Positive span of integer vectors]\label{def:span}
The positive span of a finite set of vectors $S=\{\textbf{v}_1, \textbf{v}_2, \ldots, \textbf{v}_k\}\subset\mathbb{Z}^{n}\backslash\{\textbf{0},\textbf{1}\}$ is defined as
$$pos(S):= \left\{\lambda_{0}\textbf{1}+\sum_{j=1}^{k}\lambda_{j}(\textbf{v}_{j}-\textbf{1}):\text{$\lambda_{0}\in\mathbb{R}_{>0}$, $\lambda_{j}\in\mathbb{R}_{\geq0}$ for $j=1,2,\ldots,k$}\right\}.$$
We say $S$ is positively independent if $\textbf{v}_{i}\notin pos(S\backslash{\{\textbf{v}_{i}\}})$ for $i=1$, $2$, $\ldots$, $k$.
\end{definition}
A simple consequence is that such a positive span is a convex cone. Thus, all the elements in the positive span of positive vectors are positive.
\begin{definition}[Integer frame]\label{def:frame}
Let $C$ be a convex cone in $\mathbb{R}^{n}$. A finite set $\mathcal{F}\subset\mathbb{Z}^{n}\backslash\{\textbf{0},\textbf{1}\}$ is an integer frame of $C$ if it is a positively independent set whose positive span is $C$. 
\end{definition}

Regis \cite{regis2016properties} proposed several algorithms for determining positively independent sets and positive spanning sets in the real sense. Now, we transfer them to the integer setting as Algorithm \ref{alg:pn-judge} and Algorithm \ref{alg:frame}.

\renewcommand{\algorithmiccomment}[1]{\bgroup\hfill//~#1\egroup}

\begin{algorithm} 
\caption{Positive-Negative Judge}\label{alg:pn-judge}
\begin{algorithmic}[1]
    \REQUIRE Given $S=\{\textbf{v}_1, \textbf{v}_2, \ldots, \textbf{v}_k\}\subset\mathbb{Z}^{n}\backslash\{\textbf{0},\textbf{1}\}$ and a vector $\textbf{u}\in\mathbb{Z}^{n}\backslash\{\textbf{0},\textbf{1}\}$.
    \ENSURE If $\textbf{u}$ belongs to $pos(S)$.
    \RETURN If the system $x_{0} \textbf{1} + \sum_{j=1}^{k} x_{j} \textbf{v}_{j} = \textbf{u}$ has a solution in $\mathbb{R}_{> 0}\times \mathbb{R}_{\geq 0}^{k}$.
    \COMMENT{This is solved using the CBC (Coin-or branch and cut) solver \cite{saltzman2002coin}, a mixed integer linear programming solver that handles linear constraints over positive real variables efficiently.}

\end{algorithmic} 
\end{algorithm}

\begin{algorithm} 
\caption{Frame Finder}\label{alg:frame}
\begin{algorithmic}[1]
    \REQUIRE Given $S=\{\textbf{v}_1, \textbf{v}_2, \ldots, \textbf{v}_k\}\subset\mathbb{Z}^{n}\backslash\{\textbf{0},\textbf{1}\}$.
    \ENSURE A subset of $S$ that is an integer frame of $pos(S)$.
    \STATE $\mathcal{F} \gets S$
    \FOR{$\textbf{v}\in S$}
    \IF{$\textbf{v} \in pos(\mathcal{F}\backslash\{\textbf{v}\})$} 
    \STATE $\mathcal{F} \gets \mathcal{F} \backslash \{ \textbf{v} \}$ \COMMENT{Implemented by Algorithm \ref{alg:pn-judge}}
    \ENDIF
    \ENDFOR
    \RETURN $\mathcal{F}$.
\end{algorithmic} 
\end{algorithm}


\section{Methods}\label{section:methods}

Roughly speaking, our method first express the given $q$-hypergeometric identities using the square bracket notation by Definition \ref{def:sym} and Lemma \ref{lem:extent def}. Then we ensure that each term inside the brackets is positive using Lemma \ref{lem:2inf}. By grouping these bracketed expressions following the rules in Lemma \ref{lem:qbinom}, we construct $q$-binomial coefficients, which lead to the final identities. To be more specific, We illustrate these processes using a simple example.

\begin{example}\label{ex:simple}
For $|c|<|ab|$, the $q$-Gauss sum states:
\begin{equation}
\sum_{n\geq 0}\frac{(a,b)_{n}(c/ab)^{n}}{(q,c)_{n}}=\frac{(c/a,c/b)_{\infty}}{(c,c/ab)_{\infty}}.    
\end{equation}
In order to put them in square brackets, we make some straightforward substitutions $a=q^{A}$, $b=q^{B}$ and $c=q^{C}$ where $A$, $B$ and $C$ are all integers such that $C-A-B > 0$. Then 
\begin{equation}\label{eq:exwithn}
\sum_{n\geq 0}q^{n(C-A-B)}\frac{[A,B]_{n}}{[1,C]_{n}}=\frac{[C-A,C-B]}{[C,C-A-B]}.
\end{equation}
Recalling our procedure from the beginning of this section, we must eliminate the subscripts $n$ to ensure that the terms inside square brackets are positive, as per Lemma \ref{lem:2inf}. However,in this case, whether these terms are positive or negative is uncertain. A reliable approach to find as many identities as possible is to consider any possible constraint conditions on these terms. We will explore this further in Subsection \ref{subsec:ccg}. Here, instead, we assume that $A$ and $B$ are negative, while $C$ is positive. Thus,
\begin{equation}\label{eq:ex1}
\sum_{n\geq 0}q^{n(C + n - 1)}\frac{[-A-n+1,-B-n+1,n+1,C+n]}{[-A+1,-B+1,1,C]}=\frac{[C-A,C-B]}{[C,C-A-B]}.
\end{equation}
This is exactly the form of Eq. \eqref{eq:q-hyper}. Next, by combining these bracketed expressions, we transform this identity into a binomial form. It is worth noting that such combinations are not unique, and we will discuss more practical combination methods in Subsection \ref{subsec:qbg}. Here, we only demonstrate one possible combination. Following \eqref{eq:ex1},
\begin{eqnarray}
&\Rightarrow&\sum_{n\geq 0}q^{n(C + n - 1)}\frac{[-A-n+1,n+1]}{[-A+1,1]}\frac{[-B-n+1,C+n]}{[-B+1,C]}\nonumber \\
     &&\hspace*{7cm} =\frac{[C-A,C-B]}{[C,C-A-B]}\\
\label{eqs:0}
&\Rightarrow&\sum_{n\geq 0}q^{n(C + n - 1)}{\binom{-A}{n}}_{q}\frac{[-B-n+1,C+n]}{[-B+C,1]}\frac{[-B+C,1]}{[-B+1,C]}\nonumber \\
     &&\hspace*{7cm} =\frac{[C-A,C-B]}{[C,C-A-B]}\\
\label{eqs:1}
&\Rightarrow&\sum_{n\geq 0}q^{n(C + n - 1)}{\binom{-A}{n}}_{q}{\binom{-B+C-1}{-B-n}}_{q}=\frac{[C-A,-B+1]}{[C-A-B,1]}\\ 
&\Rightarrow&\sum_{n\geq 0}q^{n(C + n - 1)}{\binom{-A}{n}}_{q}{\binom{-B+C-1}{-B-n}}_{q}={\binom{-A-B+C-1}{-B}}_{q}.
\end{eqnarray}
We begin with the sum side by combining the bracketed terms associated with $n$ into $q$-binomial coefficients. The remaining bracketed terms independent of $n$ can then be extracted from the summation and moved to the product side. In this example, the right side of Eq. \eqref{eqs:1} precisely forms a $q$-binomial coefficient. Usually, similar operations to those on the sum side are also required on the product side.
In Eq. \eqref{eqs:0}, $[-B+C]$ and $[1]$ are introduced as intermediate terms in order to generate ${\binom{-B+C-1}{-B-n}}_{q}$. This is feasible because $[-B+C]$ and $[1]$ are all positive. Unfortunately, it's not always guaranteed that the intermediate terms are positive. Therefore, we must verify the positivity of intermediate terms before incorporating them into our identity.
\end{example}

\subsection{Constraint Conditions Generator}\label{subsec:ccg}
Recall the process from Eq. \eqref{eq:exwithn} to Eq. \eqref{eq:ex1} and the addition of intermediate terms. 
Each time we encounter a $[m]$, whether it originates from our original identity or is newly introduced, we must ensure that $m$ is positive to avoid scenarios where it equals zero making the identity trivial. Sometimes, this can be confirmed or contradicted based on prior knowledge. For instance, in Eq. \eqref{eq:exwithn}, $C-A>0$ holds due to the assumption that $A<0$ and $C>0$, whereas $A-C$ cannot be positive under the same assumption. However, there are other instances where this cannot be inferred from prior knowledge alone. In such cases, we must force it to be positive to proceed with our process. As for the square bracket with positive subscript in Eq. \eqref{eq:exwithn}, say $[m]_{n}$, either $m>0$ or $m<0$ and $-m-n+1>0$ provides a feasible constraint condition due to Lemma \ref{lem:2inf}. In the case of $[m]$, things get simpler as only $m>0$ provides a feasible constraint condition. 

\begin{figure}[t]
\centering
\includegraphics[width=\textwidth]{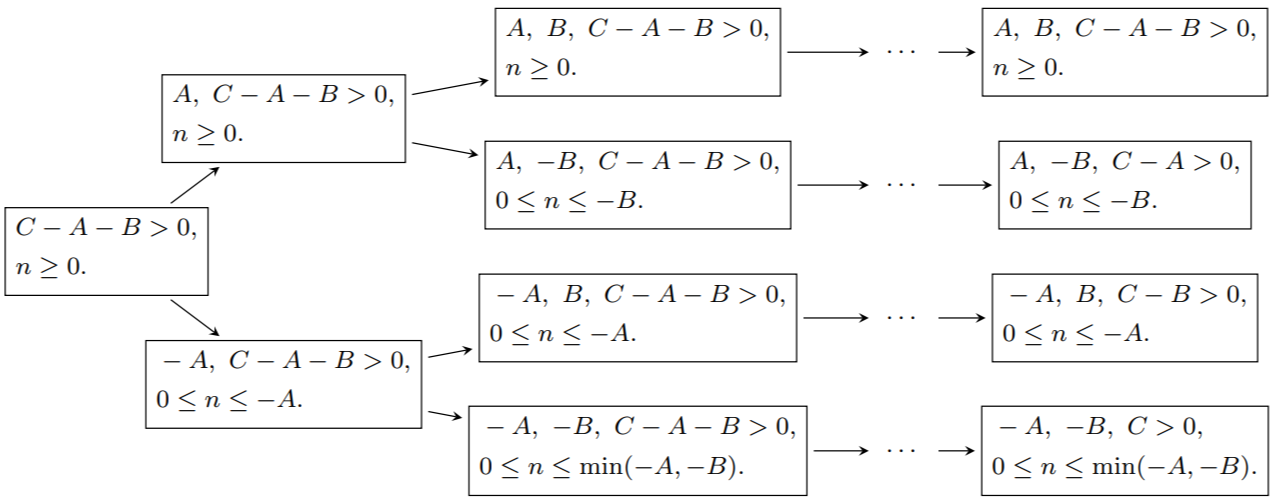}
\caption{Constraint Condition Tree}\label{fig:constraints}
\end{figure}

To help computers comprehend all the aforementioned processes, we express the variables as vectors. Continuing with Example \ref{ex:simple}, we write $m=\alpha_{0}+\alpha_{1}n+\alpha_{2}A+\alpha_{3}B+\alpha_{4}C$ as vector $\phi(m)=(\alpha_{0},\alpha_{1},\alpha_{2},\alpha_{3},\alpha_{4})$. 
Then by repeating Algorithm \ref{alg:ccg} until all brackets in Eq. \eqref{eq:exwithn} are considered or the output is empty, we will obtain all possible constraint conditions structured in a tree shape as Figure \ref{fig:constraints}. In this tree, each node represents a feasible constraint condition composed of inequalities involving integer affine expressions. Thus, each condition node corresponds to some positive integer vectors. Let $C$ denote the positive span of these vectors. Then $-C \cap C=\emptyset$. Otherwise, no variables can satisfy this constraint.

\begin{algorithm}
\caption{Constraint Conditions Generator}
\begin{algorithmic}[1]\label{alg:ccg}
\REQUIRE Given a parent condition node and a non-trivial integer affine expression $m$. 
\ENSURE All child condition nodes.
\STATE $S \gets$ the set of all corresponding positive vectors of the parent node
\STATE $\mathcal{F} \gets $ an integer frame of $pos(S)$
\STATE $\mathcal{C} \gets \emptyset$
\IF {$m$ is in a square bracket without subscript}
    \IF 
    {$\phi(m-1) \notin -pos(\mathcal{F})$}
        \STATE $\mathcal{F} \gets$ an integer frame of $pos(\mathcal{F} \cup \{\phi(m)\})$
        \STATE $\mathcal{C} \gets$ $\mathcal{C} \cup \{\text{the corresponding condition of $\mathcal{F}$}\}$
    \ENDIF
\ELSIF{$m$ is in a square bracket with positive subscript $n$}
    \STATE $\mathcal{F} \gets$ an integer frame of $pos(\mathcal{F} \cup \{\phi(n)\})$
    \IF
    {$\phi(m-1) \notin -pos(\mathcal{F})$}
        \STATE $\mathcal{F}_{1} \gets$ an integer frame of $pos(\mathcal{F} \cup \{\phi(m)\})$
        \STATE $\mathcal{C} \gets$ $\mathcal{C} \cup \{\text{the corresponding condition of $\mathcal{F}_{1}$}\}$
    \ENDIF
    \IF 
    {$\phi(m+1) \notin pos(\mathcal{F})$ and $\phi(-m-n) \notin -pos(\mathcal{F})$}
        \STATE $\mathcal{F}_{2} \gets$ an integer frame of $pos(\mathcal{F} \cup \{\phi(-m-n+1)\})$
        \STATE $\mathcal{C} \gets$ $\mathcal{C} \cup \{\text{the corresponding condition of $\mathcal{F}_{2}$}\}$
    \ENDIF
\ENDIF
\RETURN $\mathcal{C}$ (the path will be discarded if $\mathcal{C}=\emptyset$)
\end{algorithmic}
\end{algorithm}

The feasibility of all nodes can be proven by induction.
First, the statement holds trivially for the initial case where there is only one constraint. Assume that the statement holds for one node with frame $\mathcal{F}$, namely, $-pos(\mathcal{F}) \cap pos(\mathcal{F})=\emptyset$. Let $m$ be a non-trivial integer affine expression and $n$ be a positive integer. When $\phi(m-1) \notin -\text{pos}(\mathcal{F})$, it follows from Definitions~\ref{def:span} and~\ref{def:frame} that $-pos(\mathcal{F}\cup\{\phi(m)\}) \cap pos(\mathcal{F}\cup\{\phi(m)\})=\emptyset$. Otherwise, there must exists a vector $\textbf{u}\in-pos(\mathcal{F}\cup\{\phi(m)\}) \cap pos(\mathcal{F}\cup\{\phi(m)\})$. This implies that there exist positive real numbers $\alpha_{0}$ and $\beta_{0}$, and non-negative real numbers $\alpha$, $\beta$, $\alpha_{\textbf{v}}$'s and $\beta_{\textbf{v}}$'s such that 
$$\alpha_{0}\textbf{1}+\alpha(\phi(m-1))+\sum_{\textbf{v}\in\mathcal{F}}\alpha_{\textbf{v}}(\textbf{v}-\textbf{1})=-\beta_{0}\textbf{1}-\beta(\phi(m-1))-\sum_{\textbf{v}\in\mathcal{F}}\beta_{\textbf{v}}(\textbf{v}-\textbf{1}).$$
Thus, we have 
$$(\alpha_{0}+\beta_{0})\textbf{1}+(\alpha+\beta)(\phi(m-1))+\sum_{\textbf{v}\in\mathcal{F}}(\alpha_{\textbf{v}}+\beta_{\textbf{v}})(\textbf{v}-\textbf{1})=0$$
where $\alpha+\beta$ cannot be zero as $\mathcal{F}$ contains only positive integer vector. Hence, we have $\phi(m-1)\in -pos(\mathcal{F})$ which is a contradiction. Similarly, by these two definitions, we have if $\phi(m+1) \notin pos(\mathcal{F})$ and $\phi(-m-n) \notin -pos(\mathcal{F})$, then $-pos(\mathcal{F} \cup \{\phi(-m-n+1)\}) \cap pos(\mathcal{F} \cup \{\phi(-m-n+1)\})$. Therefore, due to the design of Algorithm \ref{alg:ccg}, the statement holds for all child nodes, which complete the proof.

\subsection{$q$-Binomial Coefficients Generator}\label{subsec:qbg}
To design an algorithm that transforms all square brackets to $q$-binomial coefficients and always stops after a finite number of steps, we eliminate the variables in the square brackets one by one as demonstrated in Example \ref{ex:simple}. Specifically, each time we will only focus on a target variable $x$, and then complete the following transformation.
\begin{equation*}
    \frac{[a_{1}+m_{1}x,a_{2}+m_{2}x,\ldots,a_{r}+m_{r}x]}{[b_{1}+n_{1}x,b_{2}+n_{2}x,\ldots,b_{s}+n_{s}x]} \rightarrow \frac{[c_{1},c_{2},\ldots,c_{t}]}{[d_{1},d_{2},\ldots,d_{u}]}\prod_{k}{\binom{f_{k}(x)}{g_{k}(x)}_{q}}
\end{equation*}
or abbreviated as
\begin{equation}\label{tr:total}
    \frac{\prod_{a+mx\in \mathcal{A}}[a+mx]}{\prod_{b+nx\in \mathcal{B}}[b+nx]} \rightarrow \frac{\prod_{c\in \mathcal{C}}[c]}{\prod_{d\in \mathcal{D}}[d]}\prod_{(f(x),g(x))\in\mathcal{G}}\binom{f(x)}{g(x)}_{q}
\end{equation}
where $m$'s and $n$'s are positive integers, $(a+mx)$'s, $(b+nx)$'s, $c$'s, $d$'s, $f(x)$'s and $g(x)$'s are positive integer affine expression. Moreover, the coefficients of $x$ in $a$'s, $b$'s, $c$'s and $d$'s are all zero. By Lemma \ref{lem:qbinom}, this transformation succeeds only if the sum of $m$'s equals the sum of $n$'s, which is equivalent to requiring that the $q$-hypergeometric sum is balanced. For any $a+mx$ in the numerator on the LHS and $b+nx$ in the denominator on the LHS, the following transformation can reduce the number of square brackets containing $x$.
\begin{equation}\label{tr:single}
    \frac{[a+mx]}{[b+nx]} \rightarrow \binom{b+nx-1}{a+mx-1}_{q}\frac{[1]}{[(b-a)+(n-m)x+1]}.
\end{equation}
Transformation \eqref{tr:single} succeeds if the introduced term $(b-a)+(n-m)x+1$ is positive. This can be achieved if $\phi((b-a)+(n-m)x) \notin -pos(\mathcal{F})$ where $\mathcal{F}$ is the frame corresponding to the current constraint condition. Following this transformation, the constraint condition should be updated using Algorithm \ref{alg:ccg}.
Although, in theory, we can achieve Transformation \eqref{tr:total} using a finite number of Transformation \eqref{tr:single} by determining the order of $a+mx$ and $b+nx$ involved in each step, it is still extremely time-consuming for a computer to exhaust all these possible orders, especially when $r$ or $s$ is particularly large. Therefore, we propose a simpler and more efficient algorithm (see Algorithm \ref{alg:bcg}) only considering the following three combinations.
$$\frac{[a+mx]}{[b+mx]}=\binom{b+mx-1}{b-a}_{q}\frac{[1]}{[b-a+1]},$$
$$\frac{[a+mx,a'+m'x]}{[(a+a'-1)+(m+m')x]}=\binom{a+a'-2}{a+mx-1}_{q}[1],$$
$$\ [a+mx,a'-mx]=\binom{a+a'-2}{a+mx-1}_{q}[1,a+a'-1].$$
Because of this simplification, our method does not always cover all possible identities under a fixed constraint condition. However, for each repeat loop in Algorithm \ref{alg:bcg}, there is more than one combination that meets the loop condition, which results in more than one valid transformation.

\begin{algorithm}
\caption{$q$-Binomial Coefficients Generator}
\begin{algorithmic}[1]\label{alg:bcg}
\REQUIRE LHS of Eq. \eqref{tr:total} and constraint condition frame $\mathcal{F}$.
\ENSURE RHS of Eq. \eqref{tr:total}.
\STATE $\mathcal{C}\gets\emptyset$, $\mathcal{D}\gets\emptyset$, $\mathcal{G}\gets\emptyset$

\FOR{$a+mx$ in $\mathcal{A}$}
    \IF{$m=0$}
        \STATE $\mathcal{A}\gets\mathcal{A}\backslash\{a+mx\}$, $\mathcal{C}\gets\mathcal{C}\cup\{a+mx\}$
    \ENDIF
\ENDFOR
\FOR{$b+nx$ in $\mathcal{B}$}
    \IF{$n=0$}
        \STATE $\mathcal{B}\gets\mathcal{B}\backslash\{b+nx\}$, $\mathcal{D}\gets\mathcal{D}\cup\{b+nx\}$
    \ENDIF
\ENDFOR

\REPEAT
   \STATE Find $a+mx$ in $\mathcal{A}$ and $b+nx$ in $\mathcal{B}$ such that $m=n$ and $\phi(b-a) \notin -pos(\mathcal{F})$ \label{line:condition1}
   \STATE Update $\mathcal{F}$ by Algorithm \ref{alg:ccg} with $b-a+1$ as input affine expression
   \STATE $\mathcal{A}\gets\mathcal{A}\backslash\{a+mx\}$, $\mathcal{B}\gets\mathcal{B}\backslash\{b+nx\}$, $\mathcal{C}\gets\mathcal{C}\cup\{1\}$, $\mathcal{D}\gets\mathcal{D}\cup\{b-a+1\}$
   \STATE $\mathcal{G}\gets\mathcal{G}\cup\{(b+nx-1,b-a)\}$
\UNTIL{the condition in line \ref{line:condition1} cannot be satisfied}

\REPEAT
   \STATE Find $a+mx$, $a'+m'x$ in $\mathcal{A}$ and $b+nx$ in $\mathcal{B}$ such that $n=m+m'$ and $b+1=a+a'$ \label{line:condition2}
   \STATE $\mathcal{A}\gets\mathcal{A}\backslash\{a+mx,a'+m'x\}$, $\mathcal{B}\gets\mathcal{B}\backslash\{b+nx\}$, $\mathcal{C}\gets\mathcal{C}\cup\{1\}$
   \STATE $\mathcal{G}\gets\mathcal{G}\cup\{(b+nx-1,a+mx-1)\}$
\UNTIL{the condition in line \ref{line:condition2} cannot be satisfied}

\REPEAT
   \STATE Find $a+mx$ and $a'+m'x$ in $\mathcal{A}$ such that $m+m'=0$ and $\phi(a+a'-2) \notin -pos(\mathcal{F})$  \label{line:condition3}
   \STATE Update $\mathcal{F}$ by Algorithm \ref{alg:ccg} with $a+a'-1$ as input affine expression
   \STATE $\mathcal{A}\gets\mathcal{A}\backslash\{a+mx,a'+m'x\}$, $\mathcal{C}\gets\mathcal{C}\cup\{1,a+a'-1\}$
   \STATE $\mathcal{G}\gets\mathcal{G}\cup\{(a+a'-2,a+mx-1)\}$
\UNTIL{the condition in line \ref{line:condition3} cannot be satisfied}
\STATE $\mathcal{C}\gets\mathcal{C}\backslash(\mathcal{C}\cap\mathcal{D})$, $\mathcal{D}\gets\mathcal{D}\backslash(\mathcal{C}\cap\mathcal{D})$
\IF{$\mathcal{A}=\mathcal{B}=\emptyset$}
    \RETURN $\mathcal{C}$, $\mathcal{D}$ and $\mathcal{G}$
\ENDIF
\end{algorithmic}
\end{algorithm}

Our $q$-Binomial Identities Finder also provides an enhanced algorithm \texttt{BinomialIdentity\_plus} in \texttt{qfinder.py}. It allows $q$-binomial coefficients to appear in the denominator, and non-$q$-binomial coefficients to exist in the equation. The former can be achieved simply by swapping the order of the numerator and denominator in Algorithm \ref{alg:bcg}. For the latter, we exhaust all possible orders in Transformation \eqref{tr:single} and retain the parts where the target variable cannot ultimately be converted into a $q$-binomial coefficient. For example, 
$$\frac{[x+1,x+2]}{[2x+1,1]}=\frac{[x+1,x+1]}{[2x+1,1]}\frac{[x+2]}{[x+1]}=\frac{1}{1-q^{x+1}}\binom{2x}{x}_{q}.$$

\subsection{Backtracking Framework}
Recall the Algorithms \ref{alg:ccg} and \ref{alg:bcg}. They both share a common framework for traversing feasible solutions. In Algorithm \ref{alg:ccg}, we construct a tree structure: each time a new integer expression is introduced, if it meets the predefined conditions, the branch continues; otherwise, the branch is discarded. For Algorithm \ref{alg:bcg}, we aim to pair elements from sets $\mathcal{A}$ and $\mathcal{B}$ such that after three repeat loops, all elements from both sets are used. Throughout the process, we can sequentially select elements from $\mathcal{A}$ and $\mathcal{B}$, continually testing whether the required intermediate terms meet the conditions. If they do, the process continues; if not, that choice is discarded, and we backtrack to the previous step to choose a different combination. Therefore, we can use backtracking to implement Algorithm \ref{alg:ccg} and explore feasible solutions for Algorithm \ref{alg:bcg}. Backtracking is a method where solutions are constructed incrementally by making choices at each step and undoing them if they fail to satisfy the constraints. It is akin to exploring different paths and, upon encountering a dead end, backtrack to the last decision point to try a different route. Instead of detailing how backtracking is implemented for Algorithms \ref{alg:ccg} and \ref{alg:bcg}, \footnote{Details can be found in the \texttt{signs\_generator} and \texttt{findGauss} functions in \texttt{qfinder.py} of the package.} we will illustrate the backtracking framework using Algorithm \ref{alg:back}.

\begin{algorithm}
\caption{Backtracking Framework}\label{alg:back}
\begin{algorithmic}[1]
    \REQUIRE Initial state $s_{0}$ in the state space $S$.
    \ENSURE All solutions in $S$.
    \STATE $\mathcal{R}\gets \emptyset$
    \STATE \textbf{def} \texttt{BACKTRACK}($s\in S$):
        \STATE \quad\textbf{if} $s$ is a solution \textbf{then}
            \STATE \quad \quad $\mathcal{R}\gets \mathcal{R}\cup\{s\}$
            \STATE \quad \quad \textbf{reuturn}
        \STATE \quad \textbf{end if}
        \STATE \quad \textbf{for} $c \in C$ \textbf{do}
        \COMMENT{$C$ denotes the set of all potential choices generated by $s$}
            \STATE \quad \quad \textbf{if} $c$ satisfies the constrains \textbf{then}
                \STATE \quad \quad \quad $s \gets f_{c}(s)$
                \COMMENT{$f_{c}$ updates the state via the choice $c$}
                \STATE \quad \quad \quad \texttt{BACKTRACK}($s$)
                \STATE \quad \quad \quad $s \gets f_{c}^{-1}(s)$
            \STATE \quad \quad \textbf{end if}
        \STATE \quad \textbf{end for}
    \STATE \textbf{end def}
    \STATE \texttt{BACKTRACK}($s_{0}$)
    \RETURN $\mathcal{R}$
    \end{algorithmic} 
\end{algorithm}

\section{Results}\label{section:results}

By feeding different $q$-hypergeometric identities into our framework (as implemented in the Python package), we obtained a large number of $q$-binomial identities. It is worth noting that many of these identities correspond to the same underlying formula with different parameter support. This is due to the nature of our finder, which systematically explores admissible constraints configurations for a given input identity. In this section, we mainly focus on Identity \ref{id:saaschiitz}, Identity \ref{id:suranyi} and Identity \ref{id:kummer}.

\subsection{Identity \ref{id:saaschiitz} and Identity \ref{id:suranyi}}
 
In this subsection, we will demonstrate some known $q$-binomial identities including Identity \ref{id:saaschiitz} and Identity \ref{id:suranyi}.
We start with our input $q$-hypergeometric identity.
\begin{identity} [$q$-Pfaff–Saalsch\"{u}tz Sum, \cite{andrews1984identities,zeilberger1987q,goulden1985bijective2}]\label{id:pfaff}
For positive integer $N$,
\begin{equation}
{}_3\phi_2 \left( \begin{array}{c}
a, b, q^{-N} \\
c, abq^{1-N}/c
\end{array}; q, q \right)
= \frac{(c/a,c/b)_{N}}{(c,c/ab)_{N}}
\end{equation}
\end{identity}
Identity \ref{id:pfaff} has been combinatorially proven by Andrews and Bressoud \cite{andrews1984identities}, Goulden \cite{goulden1985bijective2} and Zeilberger \cite{zeilberger1987q}. Assume $a=q^{A}$, $b=q^{B}$ and $c=q^C$ for some integers $A$, $B$ and $C$. Then feeding it to our finder leads to many $q$-binomial identities under different constraint conditions. In this subsection, we focus on two parallel constraint conditions.

\subsubsection*{Condition I}
For $-A$, $B$, $C$, $N$ and $A + B - C - N + 1 > 0$, the total six outputs are as follows.







\begin{multline}\label{res:1}
\sum_{r = 0}^{\min(-A,N)} q^{(r+A)(r-N)} {\binom{B + r - 1}{B - C}}_{q} {\binom{N}{r}}_{q} {\binom{B - C - N}{-A - r}}_{q} \\= {\binom{-A + C + N - 1}{-A}}_{q} {\binom{B - 1}{A + B - C}}_{q}.
\end{multline}

\begin{multline}\label{res:2}
\sum_{r = 0}^{\min(-A,N)} q^{(r+A)(r-N)} {\binom{B + r - 1}{B - C}}_{q} {\binom{-A}{r}}_{q} {\binom{A + B - C}{N - r}}_{q} \\= {\binom{-A + C + N - 1}{N}}_{q} {\binom{B - 1}{C + N - 1}}_{q}.
\end{multline}

\begin{multline}\label{res:3}
\sum_{r = 0}^{\min(-A,N)} q^{(r+A)(r-N)} {\binom{B + r - 1}{r}}_{q} {\binom{C + N - 1}{N - r}}_{q} {\binom{B - C - N}{-A - r}}_{q}  \\= {\binom{-A + C + N - 1}{N}}_{q} {\binom{B - C}{-A}}_{q}.
\end{multline}

\begin{multline}\label{res:4}
\sum_{r = 0}^{\min(-A,N)} q^{(r+A)(r-N)} {\binom{B + r - 1}{r}}_{q} {\binom{-A + C - 1}{-A - r}}_{q} {\binom{A + B - C}{N - r}}_{q} \\= {\binom{-A + C + N - 1}{-A}}_{q} {\binom{B - C}{N}}_{q}.
\end{multline}

\begin{multline}\label{res:5}
\sum_{r = 0}^{\min(-A,N)} q^{(r+A)(r-N)} {\binom{B + r - 1}{-A + C + N - 1}}_{q} {\binom{C + N - 1}{N - r}}_{q} {\binom{-A}{r}}_{q} \\= {\binom{B - C}{N}}_{q} {\binom{B - 1}{A + B - C}}_{q}.
\end{multline}

\begin{multline}\label{res:6}
\sum_{r = 0}^{\min(-A,N)} q^{(r+A)(r-N)} {\binom{B + r - 1}{-A + C + N - 1}}_{q} {\binom{-A + C - 1}{-A - r}}_{q} {\binom{N}{r}}_{q} \\= {\binom{B - 1}{C + N - 1}}_{q} {\binom{B - C}{-A}}_{q}.
\end{multline}

By Eq. \eqref{res:1}-\eqref{res:6}, we note that under this constraints condition, the fact that $-A$ and $N$ are interchangeable in Identity \ref{id:pfaff} is inherited. By making simple variable substitutions of $(A, B, C, D)$, we can obtain some classic $q$-binomial identities which are listed in Table \ref{tab:condition1}. 


\renewcommand{\arraystretch}{2}
\begin{table}[t]
\centering
\resizebox{\textwidth}{!}{
\begin{tabular}{|c|c|c|c|}
\hline
Identity                       & $(A,B,C,N)$          & $q$-Binomial Identitity & Reference \\ \hline
\eqref{res:1} & $(-n,3n+1,n+1,n)$    &$\sum_{k \ge 0} q^{(k-n)^2} {\binom{n}{k}}_{q} {\binom{n}{k}}_{q} {\binom{3n+k}{2n}}_{q}  = {\binom{3n}{n}}_{q} {\binom{3n}{n}}_{q}$&           \cite{gould1972combinatorial}\\ \hline
\eqref{res:2} & $(-r,t+1,t-m+1,s)$   &$\sum_{j \ge 0} q^{(j-s)(j-r)} {\binom{r}{j}}_{q} {\binom{m-r}{s-j}}_{q} {\binom{t+j}{m}}_{q}  = {\binom{t}{m-s}}_{q} {\binom{t-m+s+r}{s}}_{q} $&           \cite{Takacs1973On}\\ \hline
\eqref{res:3} & $(-a,x+y+1,-a+y+1,b)$   &Identity \ref{id:saaschiitz}&           \cite{gould1972new}\\ \hline
\eqref{res:4} & $(-a,x+y+1,y-a+1,b)$ &$\sum_{k \ge 0} q^{(a-k)(b-k)} {\binom{x+y+k}{k}}_{q} {\binom{y}{a-k}}_{q} {\binom{x}{b-k}}_{q}  = {\binom{x+a}{b}}_{q} {\binom{y+b}{a}}_{q} $&           \cite{Stanley1970Ordered}\\ \hline
\eqref{res:5} & $(m-M,m+n+1,m+1,N)$  &$\sum_{r \ge 0} q^{(N-r)(M-r-m)} {\binom{M - m}{r}}_{q} {\binom{N + m}{m + r}}_{q} {\binom{m + n + r}{M + N}}_{q}  = {\binom{m+n}{M}} {\binom{n}{N}}_{q} $&           \cite{andrews1998theory}\\ \hline
\eqref{res:6} & $(-r,x+n+r+1,1,n)$   &$\sum_{k \ge 0} q^{(k-n)(k-r)} {\binom{n}{k}}_{q} {\binom{r}{k}}_{q} {\binom{x+n+r+k}{n+r}}_{q} = {\binom{x+n+r}{n}}_{q} {\binom{x+n+r}{r}}_{q}$&           \cite{gould1972combinatorial}\\ \hline
\eqref{res:6} & $(-c-d,a+1,1-d,b)$   &$\sum_{k \ge 0}  q^{(b-k)(c+d-k)} {\binom{b}{k}}_{q} {\binom{c}{k-d}}_{q} {\binom{a+k}{b+c}}_{q}  = {\binom{a}{b-d}}_{q} {\binom{a+d}{c+d}}_{q} $&           \cite{Bizley1970AGeneralization}\\ \hline
\eqref{res:6} & $(-d,a+1,1+c-d,b)$   &$\sum_{k \ge 0} q^{(b-k)(d-k)} {\binom{b}{k}}_{q} {\binom{c}{d-k}}_{q} {\binom{a+k}{b+c}}_{q}  = {\binom{a}{b+c-d}}_{q} {\binom{a-c+d}{d}}_{q} $&           \cite{Bizley1970AGeneralization}\\ \hline
\end{tabular}%
}
\caption{Some Classic $q$-Binomial Identities}\label{tab:condition1}
\end{table}
\renewcommand{\arraystretch}{1}

\subsubsection*{Condition II}
For $-A$, $-B$, $C$, $N$ and $-A - B + C + N - 1 > 0$, there are six outputs. We only list one of them as follows since $-A$, $-B$ and $N$ are interchangeable in Identity \ref{id:pfaff} and this condition.

\begin{multline}\label{res:9}
\sum_{r = 0}^{L} q^{r(C+r-1)} {\binom{-A-B+C+N-r-1}{N-r}}_{q} {\binom{-B+C-1}{-B-r}}_{q} {\binom{-A}{r}}_{q} \\= {\binom{-A+C+N-1}{N}}_{q} {\binom{-B+C+N-1}{-B}}_{q}
\end{multline}
where $L = \min(-A,-B,N,-A - B + C + N - 1)$. Letting $(A,B,C,N)=(-k,-n,1,x)$ in Eq. \eqref{res:9} leads to Identity \ref{id:suranyi}.

\subsection{Identity \ref{id:kummer} and Its Proof}
Identity \ref{id:kummer} is derived by Bailey–Daum $q$-Kummer sum, which is stated as follows.
\begin{identity}[Bailey–Daum $q$-Kummer Sum, \cite{bailey1941note, daum1942basic,gasper2011basic}]\label{id:bailey}
For positive integer $N$,
\begin{equation}
{}_2\phi_1 \left( \begin{array}{c}
a, b \\
aq/b
\end{array}; q, -q/b \right)
= \frac{(-q;q)_{\infty}(aq,aq^{2}/b^{2};q^{2})_{\infty}}{(-q/b,aq/b;q)_{\infty}}
\end{equation}
\end{identity}
Our finder only accepts identity of the form given in Eq. \eqref{eq:exwithn}. Thus, we first set $a=q^{2A}$ and $b=q^{B}$ for some integers $A$ and $B$. Then by Lemma \ref{lem:extent def}, Identity \ref{id:bailey} is equivalent to the following identity.
\begin{equation}\label{eq:input}
    \sum_{r \geq 0} (-1)^{r}q^{r(-B+1)}\frac{[2A,B]_{r}}{[1,2A-B+1]_{r}}=\frac{[2A,-B+1][1,A-B+1;2]}{[1,2A-B+1][A,-B+1;2]}
\end{equation}
Applying \texttt{BinomialIdentity\_plus} to Eq. \eqref{eq:input} outputs Identity \ref{id:kummer}. In the algorithm, the square brackets $[m;k]$ will be grouped according to $k$ and then be processed in parallel. Moreover, following the procedure of our algorithm leads to the following proof of Identity \ref{id:kummer}.
\begin{proof}
For $A>0$ and $B<0$, Eq. \eqref{eq:input} implies
\begin{eqnarray*}
    && \sum_{r = 0}^{-B} q^{\binom{r+1}{2}}\frac{[2A,-B-r+1,2A-B+r+1,r+1]}{[2A+r,-B+1,2A-B+1,1]} \\
     &&\hspace*{5cm} =\frac{[2A,-B+1][1,A-B+1;2]}{[1,2A-B+1][A,-B+1;2]}\\
    &\Rightarrow& \sum_{r = 0}^{-B} q^{\binom{r+1}{2}} {\binom{2A+r-1}{r}}_{q}{\binom{2A-2B}{-B-r}}_{q} \frac{[2A-2B+1,1]}{[2A-B+1,-B+1]} \\
     &&\hspace*{5cm} =\frac{[2A,-B+1]}{[1,2A-B+1]}{\binom{A-B}{-B}}_{q^{2}}^{-1}\\
    &\Rightarrow& \sum_{r = 0}^{-B} q^{\binom{r+1}{2}} {\binom{2A+r-1}{r}}_{q}{\binom{2A-2B}{-B-r}}_{q}{\binom{A-B}{-B}}_{q^{2}} \\
     &&\hspace*{5cm} = \frac{[2A,-B+1,-B+1]}{[2A-2B+1,1,1]}\\
    &\Rightarrow& \sum_{r = 0}^{-B} q^{\binom{r+1}{2}} {\binom{2A+r-1}{r}}_{q}{\binom{2A-2B}{-B-r}}_{q}{\binom{A-B}{-B}}_{q^{2}} \\
     &&\hspace*{5cm} = {\binom{2A-2B}{-2B+1}}_{q} \frac{[-B+1,-B+1]}{[-2B+2,1]}\\
    &\Rightarrow& \sum_{r = 0}^{-B} q^{\binom{r+1}{2}} {\binom{2A+r-1}{r}}_{q}{\binom{2A-2B}{-B-r}}_{q}{\binom{A-B}{-B}}_{q^{2}} \\
     &&\hspace*{5cm} = {\binom{2A-2B}{-2B+1}}_{q}{\binom{-2B+1}{-B}}_{q} \frac{[-B+1]}{[-B+2]}.
\end{eqnarray*}
Our finder will continue to transform $\frac{[-B+1]}{[-B+2]}$ to $(1-q){\binom{-B+1}{1}}_{q}$, but we simply write it as $1-q^{-B+1}$ to obtain the following result.
\begin{multline}\label{eq:final}
    \sum_{r = 0}^{-B} q^{\binom{r+1}{2}} {\binom{2A+r-1}{r}}_{q}{\binom{2A-2B}{-B-r}}_{q}{\binom{A-B}{-B}}_{q^{2}} \\= \left(1-q^{-B+1}\right){\binom{2A-2B}{-2B+1}}_{q}{\binom{-2B+1}{-B}}_{q}
\end{multline}
Thus, setting $A=m$ and $B=-n$ in Eq. \eqref{eq:final} leads to Identity \ref{id:kummer}. 
\end{proof}


\section{Future Works}\label{section:further works}
Only a few identities are demonstrated in this paper. Our $q$-binomial identities finder remains hungry, waiting to be fed with more $q$-hypergeometric identities including multiple-term summation formulas \cite{Chen_Hou_Mu_2008}. Moreover, there exist many other variant generalizations of Gaussian binomial coefficients, including Fibonomial coefficients \cite{shannon2020gaussian, nived2023combinatorial}, rising binomial coefficients - type 2 \cite{shannon2020saalschutz}, Gaussian $q$-binomial coefficients with two additional parameters \cite{chu2020quadratic}, a generalization of binomial coefficients, replacing the natural numbers by an arbitrary sequence \cite{fontene1915generalisation, konvalina2000unified, he2022some}, and  the second quantized binomial coefficients \cite{kupershmidt2010mathematics, jacob2012quantum} for quantum scenario. To put them in the form of square brackets and modify our finder accordingly can also lead to related identities.

Another future direction is the opposite direction. Andrews \cite{andrews1974applications} outlined a method for translating binomial coefficient identities into hypergeometric series identities. The paper also mentioned the potential for programming this method, but it was not fully implemented. This leaves an interesting direction for future work: developing a Python implementation to translate $q$-binomial coefficient identities into $q$-hypergeometric series identities.

Although our finder is capable of generating $q$-binomial identities and providing proofs derived from $q$-hypergeometric identities, it is insufficient in combinatorial interpretation and further application. Nowadays, the $q$-binomial identities have been followed by many applications across various fields. The generating function for sets of pairs of partitions that have an ordering relation between parts of the two pairs were studied in the form of $q$-binomial coefficients in \cite{burge1993restricted}, which led to an infinite number of identities between a single and a multiple summation. Following this study, Foda et al. \cite{foda1998burge} obtained an infinite tree of $q$-polynomial identities for the Virasoro characters, by using the $q$-Saalsch{\"u}tz derived Burge transform, along with a combinatorial identity between partition pairs. Investigating the properties of $q$-binomial identities has also led to interesting results on supercongruences for $q$-analogs of the integer sequences, such as Ap{\'e}ry numbers \cite{straub2019supercongruences}. Additionally, by using Identity \ref{id:saaschiitz} and the Chinese remainder theorem for coprime polynomials, Wei et al. \cite{wei2021q} derived a series of q-supercongruences related to the third power of cyclotomic polynomials. Moreover, binomial identities were even applied to the design of privacy schemes and the assessment of their privacy performance \cite{bewong2019privacy}. Thus, newly discovered identities should not be seen as the endpoint of research, but rather as the starting point for further exploration and discovery.

\section*{Acknowledgments}
We are grateful to George Andrews for his encouraging feedback on our work. His recognition of the potential for further research in this area is much appreciated, and his earlier contributions to the field have been an important reference for our study. We also thank Iskander Aliev for some helpful comments on integer positive space, and Dawei Niu for some helpful comments on $q$-hypergeometric identities. 
This work was supported by the National Natural Science Foundation
of China (Grant No. 12301004).

\end{document}